\documentclass [12pt,a4paper,reqno]{amsart}
\usepackage{amsmath,amsthm}
\usepackage{amsfonts}
\usepackage{amssymb}
\allowdisplaybreaks

\newcommand{\be}{\begin{equation}}
\newcommand{\ef}{\end{equation}}
\chardef\bslash=`\\ 

\newtheorem*{thm*}{Theorem}

\theoremstyle{definition}

\newtheorem*{remark*}{Remarks}
\newtheorem*{defn*}{Definition}

\theoremstyle{remark}

\newcommand{\wt}{\widetilde}
\newcommand{\wh}{\widehat}
\newcommand{\fc}{\frac}

\newcommand{\iy}{\infty}
 \renewcommand{\sectionmark}[1]{}
\renewcommand{\Im}{\operatorname{Im}}
\renewcommand{\Re}{\operatorname{Re}}

\newcommand{\qc} {quasiconformal}

\newcommand{\ve}{\varepsilon}

 \usepackage{amsfonts}
\newcommand{\field}[1]{\mathbb{#1}}

\newcommand{\D}{\field{D}}
\newcommand{\om}{\omega}
\newcommand{\z}{\zeta}
\newcommand{\ov}{\overline}
\newcommand{\vp}{\varphi}
\newcommand{\hC}{\wh{\field{C}}}
\newcommand{\C}{\field{C}}
\newcommand{\R}{\field{R}}

\newcommand{\B}{\mathbf{B}}
\newcommand{\T}{\mathbf{T}}
\newcommand{\Belt}{\operatorname{Belt}}

\newcommand{\Teich}{\operatorname{Teich}}

\newcommand{\vk} {\varkappa}
\newcommand{\x} {\mathbf x}
\renewcommand{\a} {\alpha}

\begin{document}
\title {Quasiconformal deformations of nonvanishing $H^p$ functions and
the Hummel-Scheinberg-Zalcman conjecture}
\author{Samuel L. Krushkal}

\begin{abstract} Recently the author proved that the 1977 Hummel-Scheinberg-Zalcman conjecture 
on coefficients of nonvanishing $H^p$ functions is true for all $p = 2m, \ m \in \mathbb N$, i.e., for
the Hilbertian Hardy spaces $H^{2m}$. As a consequence, this also implies a proof of the Krzyz conjecture
for bounded nonvanishing functions which originated this direction.

In the present paper, we solve the problem for all spaces $H^p$ with $p \ge 2$.
\end{abstract}

\date{\today\hskip4mm({QcHSZ22.tex})}

\maketitle

\bigskip

{\small {\textbf {2020 Mathematics Subject Classification:} Primary:
30C50, 30C55, 30H05; Secondary 30F60}

\medskip


\bigskip

\markboth{Samuel L. Krushkal}{ Quasiconformal deformations and the Hummel-Scheinberg-Zalcman conjecture} \pagestyle{headings}

\bigskip
\centerline{\bf 1. INTRODUCTORY REMARKS AND MAIN RESULT}

\bigskip
The Hummel-Scheinberg-Zalcman conjecture \cite{HSZ} generalizes the
Krzyz conjecture for nonvanishing $H^\iy$ functions
to the Hardy spaces $H^p$ of holomorphic functions $f(z)$ on the unit
disk $\D = \{z : |z| < 1\}$ with norm
$$
\|f\|_p = \sup_{r<1} \Bigl(\fc{1}{2 \pi} \ \int\limits_{- \pi}^\pi |f(r e ^{i \theta})|^p d \theta \Bigr)^{1/p}.
$$
It was posed in 1977 and states that Taylor's coefficients of
nonvanishing holomorphic functions $f(z) = \sum\limits_0^\infty c_n z^n$
from $H^p, \ p > 1$,  with $\|f\|_p \le 1$ are sharply estimated by
 \be\label{1}
|c_n| \le (2/e)^{1 - 1/p},
\end{equation}
and this bound is realized only by the functions $\epsilon_2
\kappa_{n,p}(\epsilon_1 z)$, where $|\epsilon_1| = |\epsilon_2| = 1$ and
 \be\label{2}
\kappa_{n,p}(z) = \Bigl[\frac{(1 + z^n)^2}{2}\Bigr]^{1/p} \ \Bigl[\exp
\frac{z^n - 1}{z^n + 1}\Bigr]^{1-1/p}.
\end{equation}

This eminent conjecture has been investigated by many authors.
A long time, the only known results here are that the conjecture is true for $n = 1$ (Brown) and $n = 2$ (Suffridge)
as well as some results for special subclasses of $H^p$, see \cite{Br}, \cite{BW}, \cite{Su}. Brown also showed that
(1) is true for arbitrary $n \ge 2$, provided $c_m = 0$ for all $m, \ 1 \le m < (n + 1)/2$, and Suffridge also
estimated sharply the coefficient $c_1$ for $0 < p \le 1$.

Recently the author proved this conjecture for Hilbertian Hardy spaces  $H^{2 m}$ with $m \in \mathbb N$ in \cite{Kr5}, applying a new approach to the coefficient problems for holomorphic functions on the disk. This approach was presented
in \cite {Kr3}, \cite{Kr4} and involves some fundamental results of Teichm\"{u}ller space theory, especially the Bers isomorphism theorem for Teichm\"{u}ller spaces of punctured Riemann surfaces.

In the limit as $p = 2m \to \iy$, one obtains as a consequence of (1) and (2) that the coefficients of nonvanishing
$H^\iy$ functions with norm $\|f\|_\iy \le 1$ are sharply estimated by $|c_n| \le 2/e$,  which proves the initial Krzyz conjecture.

In fact, the last estimate holds for some broader class containing also unbounded nonvanishning functions, see \cite{Kr5}.

\bigskip
The aim of the present paper, which continues \cite{Kr5}, is to prove the Hummel-Scheinberg-Zalcman conjecture for all spaces $H^p$ with $p \ge 2$. The main result states:

\bigskip\noindent
{\bf Theorem 1}. {\it The estimate (1) is valid for all spaces $H^p$ with $p \ge 2$; that is, the coefficients of any
nonvanishing function $f \in \ H^p, \ p \ge 2$, with $\|f\|_p \le 1$ satisfy
$|c_n| \le (2/e)^{1 - 1/p}$ for any $n > 1$;  the
equality in (1) is realized only on the function $f(z) = \kappa_{n, p}(z)$ and its compositions with pre and post
rotations about the origin. }

\bigskip
The proof of this general theorem is based on the same ideas as its special case $p = 2 m, \ m \in \mathbb N$,
in \cite{Kr5}. A new essential step is to establish that in the case of nonvanishing $H^p$ functions the needed quasiconformal deformations exist for all $p > 1$.

Another essential step in the proof of Theorem 1 relies on the fact that for any $n \ge 2$ the coefficients of the
function $\kappa_{1,p}$ satisfy
$$
|c_n(\kappa_{1,p})| = |c_n^0| < |c_1^0|,
$$
which follows, for example, from Parseval's equality. But this remains unknown for $f \in H^p$ with $1 < p < 2$.

The last section of the paper contains some remarks concerning the possible generalizations of this problem.

\bigskip\bigskip
\centerline{\bf 2. QUASICONFORMAL DEFORMATIONS OF $H^p$ FUNCTIONS}

\bigskip
We start with establishing the existence of some special quasiconformal deformations of nonvanishning Hardy and
Bergman functions.

The general result established in \cite{Kr2} for the generic $H^p$ functions with $p = 2 m, \ m \in \mathbb N$, states the follows.

Consider the functions $f(z) \in H^{2m} \cap L_\iy(\D)$, with
$$
\sup_\D |f(z)| = M  > \|f\|_{2m}.
$$
Let $E$ be a ring domain bounded by a closed curve $L \subset \D$ containing inside the origin and by the unit circle $S^1 = \partial \D$.
Let, in addition,
$$
{\bold d}^0 = (0, 1, 0, ... \ , 0) =: (d_k^0) \in \mathbb R^{n+1},
$$
and $|\mathbf x|$ denote the Euclidean norm in $\mathbb R^l$.

\bigskip\noindent
{\bf Proposition 1}. \cite{Kr2} {\it For any holomorphic function $f(z) =
\sum\limits_{k=j}^\infty c_k^0 z^k \in L_{2m}(E) \cap L_\infty(E)$
(with $c_j^0 \ne 0, \ 0 \le j < n$ and $m \in \mathbb N$), which is
not a polynomial of degree $n_1 \le n$, there exists a
positive number $\ve_0$ such that for every point
$$
{\mathbf d}^\prime = (d'_{j+1}, \dots, \ , d'_n) \in \C^{n-j}
$$
and every $a \in \mathbb R$ satisfying the inequalities
$$
|{\mathbf  d}^\prime| \le \ve, \ \ |a| \le \ve,
\ \ \ve < \ve_0,
$$
there exists a quasiconformal automorphism $h$ of the
complex plane $\hC$, which is conformal in the disk
$$
D_0 = \{w : \ |w - c_0^0| < \sup_\D |f_0(z)| + |c_0^0| +1\}
$$
(hence also outside of $E$) and satisfies the conditions:

(i) $h^{(k)}(c_0^0) = k! d_k = k! (d_k^0 + d_k^\prime), \  k = j +
1, \dots, n$ (i.e., $d_1 = 1 + d_1^\prime$ and
$d_k = d_k^\prime$ for $k \ge 2)$;

(ii) $\|h \circ f \|_{2m}^{2m} =  \|f \|_{2m}^{2m} + a$. }

\bigskip
For any function $f \in H^p$, we have
 \be\label{3}
\|f\|_p = \sup_{r<1} \Bigl(\fc{1}{2 \pi} \ \int\limits_{- \pi}^\pi |f(r e ^{i \theta})|^p d \theta \Bigr)^{1/p} =
\Bigl(\fc{1}{2 \pi} \ \int\limits_{- \pi}^\pi |f(e ^{i \theta})|^p d \theta \Bigr)^{1/p},
\end{equation}
since the mean function
$$
\mathcal M f(r)^p = \fc{1}{2 \pi} \ \int\limits_{- \pi}^\pi |f(r e ^{i \theta})|^p d \theta
$$
is a circularly symmetric subharmonic function on $\D$,
monotone increasing with $r \to 1$.
Any such function is logarithmically convex with respect to $\log r$ and has at least one-side derivative on $[0, 1]$.

Using the appropriately thin rings $E$ adjacent to
the unit circle, one derives that for any bounded in $\D$ function $f \in H^{2m}$
there exists a $O(\ve)$-quasiconformal automorphism $h$
of $\hC$ satisfying the conditions $(i)$ and distorting the $H^{2m}$-norm of $f$ by
 \be\label{4}
\|h \circ f \|_{2m} =  \|f \|_{2m} + O(\ve)
\end{equation}
(where the bound of the remainder term depends on $m$).

The proof of Proposition 1 in \cite{Kr2} shows that all its assumptions are essential;
the arguments do not extend to arbitrary $p \ge 2$ and unbounded holomorphic $L^p$ functions.

\bigskip
One of the important steps in the proof of Theorem 1 is the following weakened extension of Lemma 1 to {\bf nonvanishing} functions from the spaces $H^p$ with $p > 1$.

Fix a natural $n >  1$ and consider the collections
$\mathbf d = (d_0, d_1, \dots, d_n) \in \mathbb C^{n+1}$.

\bigskip\noindent
{\bf Proposition 2}. {\it For every bounded nonvanishing function $f(z) = c_0 + c_1 z + c_2 z^2 + \dots \in H^p, \ p > 1$, with $\|f\|_p < \iy$, which is not a polynomial of degree at most $n$, there exists $\ve_0 > 0$ such that for any point $\mathbf d \in \mathbb C^{n+1}$ with $|\mathbf d| \le \ve < \ve_0$, there is a bounded nonvanishing function $f^*(z) = c_0^* + c_1^* z + c_2^*  z^2 + \dots \in H^p$ satisfying $c_j^* = c_j + d_j$ for all $j = 0, 1, \dots, n$ and }
$$
\|f^*\|_p = \|f\|_p + O(\ve).
$$

\bigskip
Note that quasiconformal deformations of such type preserve the $L_p$-norm of holomorphic functions, but generically increase their $L_\iy$-norm (an illustrating example is given in \cite{Kr4}).

\bigskip\noindent
{\bf Proof}. For any nonvanishing function $f(z) = \sum\limits_0^\iy c_n z^n \in H^p, p > 1$, the function
$$
f_{p/2}(z) := f(z)^{p/2} = e^{(p/2) \log f(z)},
$$
with a fixed branch of the logarithmic function, also is single valued, holomorphic and zero free in the unit disk $\D$.
We take everywhere the principal branch.
Explicitly,
 \be\label{5}
f_{p/2}(z) = c_0^{p/2} \Bigl(1 + \fc{p}{2} \  \fc{c_1}{c_0} \  z + \dots \Bigr) =  c_0(f_{p/2}) + c_1(f_{p/2}) z + \dots \ ;
\end{equation}
this function belongs to the space $H^2$. In the case of nonvaninshing $f$, the correspondence
$f(z) \leftrightarrow f_{p/2}(z)$ creates a biholomorphism (one-to-one and open map) between the neighborhoods of the origin
in $\C^{n+1}$ filled by collections $\mathbf d(f) = (c_0, \dots, c_n)$ and $\mathbf d(f_{/2})p = (c_0(f_{p/2}), \dots, c_n(f_{p/2}))$.

Since the function (5) is holomorphic, one can apply to it all arguments used in \cite{Kr2} in the proof of Proposition 1 for $p = 2m$ (this proof essentially requires that $f(z)^m$ is single valued).
In view of the importance of Proposition 2, we outline the main steps of its proof; the details omitted are given in
\cite{Kr2}.

Fix $R \ge \sup_\D |f(z)| + |c_0| +1$
and take the annulus
$$
G_R = \{w : \ R < |w - c_0^0| < R + 1\}.
$$
We define  for $\rho \in L_p(B), \ p \ge 2$, the operators
$$
T\rho = - \frac{1}{\pi} \iint\limits_{G_R} \frac{\rho(\zeta) d \xi d \eta}{\z - w}, \quad
\Pi \rho = \partial_w T\rho = - \frac{1}{\pi} \iint\limits_{G_R}  \frac{\rho(\zeta) d \xi d \eta}{(\z - w)^2}
$$
(the second integral exists as a principal Cauchy value).
We seek the required quasiconformal automorphism $h = h^\mu$ of the form
 \be\label{6}
h(w) = w - \frac{1}{\pi} \iint\limits_{G_R}  \frac{\rho(\zeta) d \xi d \eta}{\zeta - w} = w + T \rho(w),
\end{equation}
with the Beltrami coefficient $\mu = \mu_h$ equal to zero outside of $G_R$, with $\|\mu \|_\iy < \kappa < 1$.
Substituting (6) into the Beltrami equation
$\partial_{\ov w} h = \mu \partial_w h$, we get
$$
\rho = \mu + \mu \Pi \mu + \mu \Pi (\mu\Pi \mu) + ... \ .
$$
This series is convergent in $L_p(E)$ for some $p > 2$, and on the basis of well-known properties of operators $T$
and $\Pi$, we have for any disk
$\D_{R^\prime} = \{w \in \C : \ |w| < R^\prime\}, \ \ 0 < R^\prime < \iy$,
that
$$
\|\rho\|_{L_p(\D_{R^\prime})},  \ \ \| \Pi \rho \|_{L_p(\D_{R^\prime})} \le  M_1(\kappa, R^\prime, p) \|\mu\|_{L_\iy(\C)};
\|h\|_{C(\D_{R^\prime})} \le  M_1(\kappa, R^\prime, p) \|\mu \|_\iy.
$$
Therefore,
$$
h(w) = w + T \mu(w) + \om(w),
$$
with $\|\om \|_{C(\D_{R^\prime})} \le M_2(\kappa, R^\prime) \|\mu \|_\iy^2$.
Using the pairing
$$
<\nu, \vp>_{G_R} = - \fc{1}{\pi} \iint\limits_{G_R} \nu(\zeta) \vp(\zeta) d \xi d \eta,
\quad \nu \in L_\iy(G_R), \ \ \vp \in L_1(G_R),
$$
one can rewrite the above representation in the form
 \be\label{7}
h(w) = w + \sum\limits_0^\iy <\mu, \vp_k>_E (w - c_0^0)^k + \om(w), \quad
\vp_k(\z) = \frac{1}{(\z - c_0^0)^{k+1}}.
\end{equation}
This equality and the condition $c_j^* = c_j + d_j$ for all $j = 0, 1, \dots, n$, provide  the first group
of equalities to determine the desired Beltrami coefficient $\mu$:
 \be\label{8}
k! d_k = <\mu, \vp_k>_{G_R} + \om^{(k)}(c_0^0) = <\mu, \vp_k>_{G_R} + O(\|\mu \|_\infty^2),
\quad k = j + 1, ... \ , n.
\end{equation}
On the other hand, (7) and the requirement of preserving $L_p$ norms give
$$
\begin{aligned}
\|h \circ f \|_p^p &= \| f + T \rho \circ f \|_p^p = \int\limits_{G_R} |f(z) + T \mu \circ f(z)|^p d E_z + O(\|\mu \|_\iy^2)  \cr
&= \int\limits_{G_R} [|f(z)|^2 +  2 \Re (\ov{f(z)} T \mu \circ f(z))
+ |T \mu \circ f(z)|^2 \bigr]^{p/2} dx dy + O(\|\mu \|_\iy^2)  \cr
&= \|f\|_p^p + \frac{p}{2 \pi} \Re \Bigl[ \iint\limits_{G_R} \mu(\z) d \xi d \eta
\int\limits_{G_R} \fc{|f(z)|^{p-2} \ \ov{f(z)}}{\z - f(z)} dx dy \Bigr] + O_p(\|\mu \|_\iy^2)
\end{aligned}
$$
(here $z = x + iy$). Now set
 \be\label{9}
\phi(\zeta) = - \fc{p}{2} \int\limits_{G_R} \fc{|f(z)|^{p-2} \ov{f(z)}}{f(z) - \zeta} dx dy;
\end{equation}
then the previous equality can be rewritten in the form
  \be\label{10}
\|h \circ f\|_p^p - \|f\|_p^p = \Re <\mu, \phi>_{G_R} + \ O_p(\|\mu \|_\iy^2).
\end{equation}

The function (9) is holomorphic in the disk $\D_R^* = \{w \in \hC : \ |w - c_0^0| > R\}$ and belongs to the space
$B_R^p$ formed in $B_p(B)$ by functions holomorphic in $D_R^*$; moreover,
$\phi(\z) \not\equiv 0$.
The latter follows from the fact that for large $|\z|$ we have
$\phi(\z) = \sum\limits_1^{\iy} b_k \z^{- k}$   with $b_2 = \fc{p}{2} \ \|f\|_p^p > 0$.

We shall need the following important lemma whose proof straightforwardly follows the corresponding lemma in \cite{Kr2}.

\bigskip\noindent
{\bf Lemma 1}. {\it Under the assumptions of the theorem, the function $\phi$ is distinct from a linear combination
of the fractions} $\vp_0, ... , \vp_l$, with $l \le n$.

\bigskip
According to Lemma 1, the series expansion of $\phi$ in $D_R^{*}$ must contain the powers $(\z - c_0^0)^{-k-1}$
with $k > n$, and therefore, the remainder
$$
\psi(\z) = \phi(\z) - \sum\limits_0^n b_k (\z - c_0^0)^{-k-1} =
\Bigl(\sum\limits_0^{j-1} + \sum\limits_s^\iy\Bigr)
b_k (\z - c_0^0)^{-k-1}, \ \ s \ge n + 1,
$$
is distinct from zero in $D_R^{*}$. Let us note also that (3) implies
$$
h \circ f(z) = c_0^* + \sum\limits_j^\iy c_k^* z^k,
$$
with
$c_0^* = c_0^0 + d'_0$ and $c_j^* = c_j^0 d'_1 \ (j \ge 1)$. Thus,
$$
|c_j^*|^2 - |c_j^0|^2 = \begin{cases} \ 2 \Re ({\ov c}_0^0 d'_0)
= 2 \Re ({\ov c}_0^0 <\mu, \vp_0>)  + O(\| \mu \|^2), &\quad j = 0, \\
2 |c_j^0|^2 \Re d'_1 = 2 |c_j^0|^2 \Re <\mu, \vp_0>  + O(\|\mu \|^2), &\quad j \ge 1,
\end{cases}
$$
and, therefore, $b_j \ne 0$.

Let us now seek the desired Beltrami coefficient $\mu$ in the form
 \be\label{11}
\mu = \xi_j {\ov \vp}_j + \sum\limits_{j+1}^n \xi_k {\ov \vp}_k + \tau \ov \psi, \ \  \mu|\C \setminus B = 0,
\end{equation}
with unknown constants
$\xi_j, \xi_{j+1}, ... \ , \xi_n, \tau$ to be determined from  equalities (8) and (10).

Substituting the expression (11) into (8) and (10) and taking into account the mutual orthogonality of $\vp_k$ on $B$,
one obtains the nonlinear equations
 \be\label{12}
\begin{aligned}
k! d_k &= \xi_k r_k^2 + O(\Vert \mu \Vert^2), \quad k = j + 1, ... \ , n, \cr
\Vert h \circ f_0 \Vert_{2m}^{2m} - \Vert f_0 \Vert_{2m}^{2m} &=
\Re <\xi_j {\bar \vp}_j + \sum\limits_{j+1}^n \xi_k {\bar \vp}_k +
\tau \bar \psi, \phi> + O(\Vert \mu \Vert^2)
\end{aligned}
\end{equation}
for determining $\xi_k$ and $\tau$. The only remaining equation is a relation for
$\Re \xi_j, \Im \xi_j, \Re \tau, \Im \tau$. To distinguish a unique solution, we add three real equations to (12).
Let us require that $\xi_j$ satisfy the equality
 \be\label{13}
<\xi_j {\ov \vp}_j + \sum\limits_{j+1}^n \xi_k {\bar \vp}_k, \sum\limits_0^n  b_k \vp_k > = 0;
\end{equation}
it is reduced to
 \be\label{14}
\xi_j b_j r_j^2 = - \sum\limits_{j+1}^n \xi_k b_k r_k^2 \quad (b_j \neq 0).
\end{equation}
Then we obtain for $\tau$ the equation
$$
\|h \circ f_0 \|_p^p - \|f \|_p^p =
\Re <\tau \ov \psi, \phi> + O(\|\mu \|_]\iy^2),
$$
which, letting $\tau$ be real, takes the form
 \be\label{15}
\|h \circ f \|_p^p - \|f\|_p^p = \tau \vk + O(\|\mu\|_\iy^2)
\end{equation}
with $\vk = \sum\limits_k r_k^2$. The summation is taken here over all $k \ne j + 1, ... \ , n$, for which $b_k \ne 0$.

Separating the real and imaginary parts in equalities (12), (14) and adding (15), we obtain $2(n - j) + 3$ real equalities, which define a nonlinear $C^1$ smooth (in fact, $\R$-analytic) map
$$
\bold y = W(\mathbf x) = W^\prime(\mathbf 0) \mathbf x + O(|\mathbf x|^2),
$$
of the points
$\mathbf x = (\Re \xi_j, \Im \xi_j, \Re \xi_{j+1}, \Im \xi_{j+1}, ... \ , \Re \xi_n, \Im \xi_n, \tau)$
in a small neighborhood $U_0$ of the origin in $\R^{2(n - j) + 3}$, taking
the values
$$
\mathbf y = (\Re d_j, \Im d_j, \Re d_{j+1}, \Im d_{j+1}, ... \ , \Re d_n, \Im d_n, \|h \circ f \|_p^p - \|f_0 \Vert_p^p)
$$
also near the origin of $\R^{2(n - j) + 3}$. Its linearization $\mathbf y = W^\prime(\mathbf 0) \mathbf x$
defines a linear map $\R^{2(n -j) + 3} \to \R^{2(n - j) + 3}$ whose Jacobian only differs from $r_j^2 r_{j+1}^2 ... r_n^2 \vk \ne 0$ by a constant factor.
Therefore, $\mathbf x \mapsto W^\prime(\mathbf 0) \bold x$ is a linear isomorphism of the space $\R^{2(n - j) + 3}$ onto itself, and one can apply to $W$ the inverse mapping theorem. The latter implies the assertion of Proposition 2.

So, for any collection $\mathbf d_\ve(f) = \mathbf d(f) + O(\ve)$ there exists an $O(\ve)$-quasiconformal homeomorphism $h^{\mu}$ of $\hC$ conformal on $f(\D)$ such that
$\mathbf d_\ve(f) = (c_0(h^{\mu} \circ f), \dots, c_j(h^{\mu} \circ f)$, and
  \be\label{16}
\|h^\mu \circ f\|_{L_p} = \|f\|_{L_p}^p.
\end{equation}

Note also that the $H^p$-norm is distorted similar to (4) via
 \be\label{17}
\|h^\mu \circ f\|_{H^2}^2 = \|f_{p/2}\|_{H^2}^2 = \|f\|_{H^p}^p + O(\ve).
\end{equation}
The relations (6) and (7) show the difference under the actions of deformations given by Propositions 1 and 2 on the Hardy
and Bergman spaces.

\bigskip\bigskip
\centerline{\bf 3. PROOF OF THEOREM 1}

\bigskip\noindent
{\bf Step 1: Underlying lemmas}.
As was mentioned in the introduction, the proof of the theorem for functions from $H^p$ with $p \ge 2$ follows the lines
of \cite{Kr5} with applying in the needed places Lemma 2.

Denote the unit ball of $H^p$ by $B_1(H^p)$ and its subset of nonvanishing functions by
$B_1^0(H^p)$. It will be convenient to regard the free coefficients
$c_0(f)$ also as elements of $B_1^0(H^p)$, which are constant on the
disk $\D$. Let
$$
\wh B_1^0(H^p) = B_1^0(H^p) \cup \{f_0\},
$$
where $f_0(z) \equiv 0$.

We shall essentially use Brown's result quoted above and present it as

\bigskip\noindent
{\bf Lemma 2}. \cite{Br} {\it For any $f(z) = c_0 + c_1 z + c_2 z^2
+ \dots \in B_1^0(H^p)$, we have
$$
|c_1| \le (2/e)^{1 - 1/p},
$$
with equality only for the rotations of function $\kappa_1(z)$ given by (2). }

\bigskip
The following important lemma concerns one of the basic intrinsic features of nonvanishing holomorphic functions
(the openness)

\bigskip\noindent
{\bf Lemma 3}. \cite{Kr5} {\it Every point $f \in B_1^0(H^p)$ has a neighborhood
$U(f, \epsilon)$ in $H^p$, which entirely belongs to $ B_1^0(H^p)$,
i.e., contains only nonvanishing $H^p$ functions on the disk $\D$.
Take the maximal balls $U(f, \epsilon)$ with such property. Then
their union
$$
\mathcal U^p = \bigcup_{f\in B_1^0(H^p)} U(f, \epsilon)
$$
is an open path-wise connective set, hence a domain, in the space
$\wh B_1^0(H^p)$. }

\bigskip
Let $\mathcal P_n$ be the linear space of polynomials of degree less than or equal to $n$, and
$\mathcal P = \bigcup_n \mathcal P_n$.

\bigskip\noindent
{\bf Lemma 4}. {\it The intersection $\mathcal U^p \bigcap \mathcal P$ is dense in $\mathcal U^p$, which means that any $f$ from the distinguished domain $\mathcal U^p$ is approximated in $H^p$ by nonvanishing polynomials. }

\bigskip
We shall use in the proof of Theorem 1 somewhat different (up to a biholomorphic homeomorphism) model of the universal  Teichm\"{u}ller space $\T$, which involves quasiconformally extendable univalent functions in the disk satisfying some non-standard prescribed normalization conditions. Their existence of such maps is ensured by the following lemma related to  solutions of the Beltrami equation
$\partial_{\ov z} w = \mu(z) \partial_z w$ on $\C$ with
coefficients $\mu$ supported in the disk $\D^*$, i.e., from the ball
$$
\Belt(\D^*)_1 = \{\mu \in L_\iy(\C): \ \mu|\D = 0, \ \|\mu\| < 1\}.
$$

\bigskip\noindent
{\bf Lemma 5}. \cite{Kr5} {\it For any Beltrami coefficient $\mu \in \Belt(\D^*)_1$ and any $\theta_0 \in [0, 2 \pi]$,
there exists a point $z_0 = e^{i \a}$ located on $\mathbb S^1$ so that
$|e^{i \theta_0} - e^{i \a}| < 1$ and such that for any $\theta$ satisfying
$|e^{i \theta} - e^{i \a}| < 1$ the equation
$\partial_{\ov z} w =  \mu(z) \partial_z w$
has a unique homeomorphic solution $w = w^\mu(z)$, which is holomorphic on the unit disk $\D$
and satisfies
 \be\label{18}
w(0) = 0, \quad w^\prime(0) = e^{i \theta}, \quad w(z_0) = z_0.
\end{equation}
This solution is holomorphic on the unit disk $\D$, and hence,  =$w^\mu(z_{*}) = \iy$ at some point $z_{*}$
with $|z_{*}| \ge 1$.  }

\bigskip\noindent
{\bf Step 2: Holomorphic embedding of nonvanishing $H^p$ functions into Teichm\"{u}ller spaces and lifting the functional $J_n(f) = c_n$}.
Denote by $\B = \B(\D)$ the space of
hyperbolically bounded holomorphic functions $\vp(z)$ (regarded as
holomorphic quadratic differentials $\vp(z) dz^2$ so that $\vp \circ h(z)h^\prime(z)^2 = \vp(z)$ for any conformal coordinate map $h$) on the unit disk, with norm
$$
\|\vp \|_\B = \sup_\D(1 - |z|^2)^2 |f(z)|.
$$
Every $\vp \in \B$ is the {\bf Schwarzian derivative}
$$
S_w(z) = \left(\frac{w^{\prime\prime}(z)}{w^\prime(z)}\right)^\prime
- \frac{1}{2} \left(\frac{w^{\prime\prime}(z)}{w^\prime(z)}\right)^2,
\quad z \in \D,
$$
of a locally
univalent function $w(z)$ in the disk $\D$ determined (up to a Moebius map of the sphere $\hC$) from the nonlinear differential equation
$$
w^{\prime\prime\prime}/w^\prime - 3 (w^{\prime\prime}/w^\prime)^2/2 = \vp,
$$
or equivalently, as the ratio $w = \eta_2/\eta_1$
of two linearly independent solutions of the linear equation $2 \eta^{\prime\prime} + \vp \eta = 0$ in $\D$.

The space $\B$ is dual to the space $A_1(\D)$ of
integrable holomorphic functions on $\D$ with $L_1$ norm.

The Schwarzians $S_w$ of functions $w$ univalent in the whole disk $\D$ and having quasiconformal extensions to $\hC$ fill a path-wise bounded domain in $\B$; this domain the most appliable model of  the {\bf universal Teichm\"{u}ller space} $\T = \Teich(\D)$ (with appropriate normalization of maps $w$).

\bigskip
Let $A_p(\D), p \ge 1$, be the Bergman spaces of holomorphic functions in $\D$ with norm
$$
\|f\|_{A_p} = \Bigl( \fc{1}{\pi} \ \iint_\D |f(z)| dx dy
\Bigr)^{1/p} \quad (z = x + i y).
$$
For each $f \in H^p$, we have
$\|f\|_{A_p}^p =  \le \fc{1}{2} \|f\|_{H^p}^p$,
which yields, since $A_p(\D) \subset A_1(\D) \subset \B$ and $\|f\|_\B \le \|f\|_{A_1(\D)}$ for $\vp \in A_1(\D)$,
that all functions $f \in H^p$ belong to the space $\B$. Therefore, these functions can be regarded as the Schwarzian derivatives of locally univalent functions in $\D$.

In particular, the functions $f$ from the ball
$$
B_\rho(H^p) = \{f \in H^p: \ \|f\| < \rho\}
$$
with radius $\rho = 1/2^{1/p}$ satisfy $\|f\|_\B < 2$, and hence
are the
Schwarzians of univalent functions in the whole disk $\D$ admitting quasiconformal extension to
the complementary disk
$$
\D^* = \{z \in \hC: \ \ |z| > 1\}.
$$
Therefore, such $f$ are the points of the universal Teichm\"{u}ller space $\T$.
This implies a {\it holomorphic embedding} $\iota$ of the ball $B_\rho(H^p)$ and of its open subset
$$
\fc{1}{2^p} \mathcal U^p = \{\fc{1}{2^p} f: \ f \in \mathcal U^p \}
$$
into the space $\T$.

\bigskip
Now consider the family $\wh S(1)$ of quasiconformally extendable to $\hC$  holomorphic univalent functions
$$
w(z) = a_1 z + a_2 z^2 + \dots, \quad z \in \D,
$$
with $|a_1| = 1$ and $w(z_0) = z_0$ for some point $z_0 \in \mathbb S^1$ (depending on $w$), completed in the topology of locally uniform convergence on $\C$.
This collection is a disjunct union
$$
\wh S(1) = \bigcup_{- \pi \le \theta < \pi} S_\theta,
$$
where $S_\theta$ consists of quasiconformally extendable univalent functions on $\D$ with expansions
$$
w(z) = e^{i \theta} z + a_2 z^2 + \dots
$$
having a fixed point $z_0 \in \mathbb S^1$ (also completed in the indicated weak topology). These collections
preserve conjugation with rotations $z \mapsto e^{i \a}z$, i.e., contain for each $w$ the rotated functions
$w_{\a,\a}(z) = e^{- i \a} w(e^{i \a} z)$. There is often enough to deal with the class $S_0$ related to $z_0 = 1$.

The assertion of Lemma 5 is also valid for the limit functions of sequences $\{w_n\}$
of functions $w_n \in \wh S(1)$ with quasiconformal extension, but in the general case
the equality $w(z_0) = z_0$ must be understand in terms of the Carath\'{e}odory prime ends.
As was indicated above, any function from $\wh S(1)$ with $\theta$ chosen following
Lemma 5 is holomorphic on the disk $\D$ (has there no pole).

This family $\wh S(1)$ is closely related to the canonical class $S$ of
univalent functions $w(z)$ on $\D$ normalized by $w(0) = 0, \
w^\prime(0) = 1$. Every $w(z) \in S$ has its representatives $w_{\tau, \theta}$ in
$\wh S(1)$ obtained by pre and post compositions of $w$ with rotations $z \mapsto e^{i \tau} z$ about
the origin, related by
 \be\label{19}
w_{\tau, \theta}(z) =  e^{- i \theta} w(e^{i \tau} z) \quad
\text{with} \ \ \tau = \arg z_0,
\end{equation}
where $z_0$ is a point of the circle $\mathbb S^1$ whose image $w(z_0) = e^{i \theta}$ is a common
point of the unit circle and the  boundary of domain $w(\D)$.

This is trivial for the identity map $w(z) \equiv z$ (then one can take $\theta = \tau = 0$).
For any another $w(z)$ the existence of such a point $z_0$ follows from the Schwarz lemma, which yields, together with the assumption $w^\prime(0) = 1$, that the image $w(\D)$ cannot lie entirely in $\D$; hence, its boundary $\partial w(\D)$ has common points with the circle $\mathbb S^1$.

This connection also implies that the functions conformal in the closed disk $\ov \D$ are dense in each class $S_\theta$.
Note also that the classes $S_\theta$ and $\wh S(1)$ are compact in the topology of locally uniform  convergence on $\D$.

\bigskip
The Schwarzian derivatives of $w$ and $w_{\tau, \theta}$ are related by
$$
S_{w_{\tau, \theta}}(z) = S_w(e^{i \tau} z) e^{2i \tau},
$$
which yields that {\it for any fixed $\theta$ the Schwarzians $S_w$ of $w \in S_\theta$ fill the same
bounded domain in the space $\B$, which models $\T$}.

In other words, the relation (9) allows us to model the universal Teichm\"{u}ller space $\mathbf T$ for any fixed $\theta$ by the Schwarzians $S_w = \vp$ of functions $w(z) = e^{i \theta} z + a_2 z^2 + \dots$ from the sets $S_\theta$.

In this case, going to the limit $\lim\limits_{t \to 0} \|S_w(t e^{i \a} z)\|_\B \to 0$ along a curve
$\{S_w(t e^{i \a} z): \ 0 \le t \le 1\}$ with fixed nonzero $\theta$ and $\a$ one attains in the space $\T$ its base point $\vp = \mathbf 0$, and the corresponding function in $ S_\theta$ is the elliptic fractional linear transformation
$$
w = \frac{e^{i \theta} z}{(1 - e^{- i \theta})z_0^{-1} z + 1}
$$
with fixed points $0$ and $z_0 = e^{i \a}$. For $\a = \theta = 0$, this is the identity map.

Note also that the relation (19) is compatible  with existence and uniqueness of appropriate conformal and  quasiconformal maps, holomorphy of their Taylor coefficients, the Teichm\"{u}ller space theory, etc.
Actually we deal with the classical model of Teichm\"{u}ller spaces via domain in the Banach spaces of Schwarzian
dervatives $S_w$ in $\mathbb D$ (or in the disk $\mathbb D^*$) of univalent holomorphic functions normalized
either by fixing three boundary points on the unit circle $S^1$ or via $w(0) = 0, \ w^\prime(0) = 1, \ w(\iy) = \iy$
(often the disk is replaced by the half-plane).

\bigskip
An equivalent model of $\T$ is obtained by applying the inverted functions
 $W(z) = 1/w(1/z)$ for $w \in S_\theta$, which form the corresponding classes $\Sigma_\theta$ of nonvanishing univalent functions on the disk $\D^*$ with expansions
$$
W(z) =  e^{- i \theta} z + b_0 + b_1 z^{-1} + b_2 z^{-2} + \dots, \quad  W(1/\a) = 1/\a,
$$
and $\wh \Sigma(1) = \bigcup_\theta \Sigma_\theta$.

Simple computations yield that the coefficients $a_n$ of $f \in S_\theta$ and the corresponding coefficients $b_j$ of $W(z) = 1/f(1/z) \in \Sigma_\theta$ are related by
$$
b_0 + e^{2i \theta} a_2 = 0, \quad b_n + \sum \limits_{j=1}^{n}
\epsilon_{n,j}  b_{n-j} a_{j+1} + \epsilon_{n+2,0} a_{n+2} = 0, \quad n = 1, 2, ... \ ,
$$
where $\epsilon_{n,j}$ are the entire powers of $e^{i \theta}$ \ ($\theta$ is fixed). This
successively implies the representations of $a_n$ by $b_j$ via
 \be\label{20}
a_n = (- 1)^{n-1} \epsilon_{n-1,0}  b_0^{n-1} - (- 1)^{n-1} (n - 2)
\epsilon_{1,n-3} b_1 b_0^{n-3} + \text{lower terms with respect to}
\ b_0.
\end{equation}

By abuse of notation, we shall denote the holomorphic embedding of $H^p$ into the space $\T$ modelled by Schwarzians in $\D^*$ by the same letter $\iota$. The image $\iota H^p$ is a non-complete linear
subspace in $\B$, and the image of the distinguished domain $\fc{1}{2^p} \mathcal U^p$ is a complex submanifold in $\T$.

Note that the coefficients $\a_n$ of Schwarzians
$$
S_w(z) = \sum_0^\infty \a_n z^n
$$
are represented as polynomials of $n + 2$ initial coefficients of $w \in S_\theta$ and, in view of (10), as polynomials of $n + 1$ initial coefficients of the corresponding $W \in \Sigma_\theta$
(provided that $\theta$ and $\a$ are given and fixed and the number $e^{i \theta}$ is considered to be a constant).

We denote these polynomials by $J_n(w)$ and $\wt J_n(W)$, respectively,
and will deal with these polynomial functionals only on the union of admissible classes
$S_\theta$ or $\Sigma_\theta$.

\bigskip\noindent
{\bf Step 3: Lifting to covering space $\T_1$ and estimating the restricted plurisubharmonic functional}.
Our next step is to lift both polynomial functionals $J_n(w)$ and
$\wt J_n(W)$ onto the Teichm\"{u}ller space $\T_1$ of the punctured disk $\D_{*} =
\D \setminus \{0\}$, which covers $\T$.

Recall that the points of $\T_1$ are the classes $[\mu]_{\T_1}$ of $\T_1$-{\bf equivalent} Beltrami
coefficients $\mu \in \Belt(\D)_1$ so that the corresponding \qc \
automorphisms $w^\mu$ of the unit disk coincide on both boundary
components (unit circle $\mathbb S^1$ and the puncture
$z = 0$) and are homotopic on $\D \setminus \{0\}$. This space also is a complex Banach
manifold.

Due to the Bers isomorphism theorem \cite{Be}, the space $\T_1$ is
biholomorphically isomorphic to the {\bf Bers fiber space}
$$
\mathcal F(\T) = \{(\phi_\T(\mu), z) \in \T \times \C: \ \mu \in
\Belt(\D)_1, \ z \in w^\mu(\D)\}
$$
over the universal space $\T$ with holomorphic projection $\pi(\psi,
z) = \psi$.
This fiber space is a bounded hyperbolic domain in $\B \times \C$
and represents the collection of domains $D_\mu = w^\mu(\D)$ as a
holomorphic family over the space $\T$.

The indicated isomorphism between $\T_1$ and $\mathcal F(\T)$ is
induced by the inclusion map \linebreak $j: \ \D_{*} \hookrightarrow
\D$ forgetting the puncture at the origin via
 \be\label{21}
\mu \mapsto (S_{w^{\mu_1}}, w^{\mu_1}(0)) \quad \text{with} \ \
\mu_1 = j_{*} \mu := (\mu \circ j_0) \ov{j_0^\prime}/j_0^\prime,
\end{equation}
where $j_0$ is the lift of $j$ to $\D$.

Now,letting
 \be\label{22}
\wh J_n(\mu) = \wt J_n(W^\mu),
\end{equation}
we lift these functionals  from the sets $S_\theta$ and
$\Sigma_\theta$ onto the ball $\Belt(\D)_1$. Then, under the
indicated $\T_1$-equivalence, i.e., by the quotient map
$$
\phi_{\T_1}: \ \Belt(\D)_1 \to \T_1, \quad \mu \to [\mu]_{\T_1},
$$
the functional $\wt J_n(W^\mu)$ is pushed down to a bounded holomorphic
functional $\mathcal J_n$ on the space $\T_1$ with the same range domain.

Equivalently, one can apply the quotient map $\Belt(\D)_1 \to \T$
(i.e., $\T$-equivalence) and compose  the descended functional on
$\T$ with the natural holomorphic map $\iota_1: \ \T_1 \to \T$
generated by the inclusion $\D_{*} \hookrightarrow \D$ forgetting
the puncture. Note that since the coefficients $b_0, \ b_1, \dots$
of $W^\mu \in \Sigma_\theta$   are uniquely determined by its
Schwarzian $S_{W^\mu}$, the values of $\mathcal J_n$ in the points
$X_1, \ X_2 \in \T_1$ with $\iota_1(X_1) = \iota_1(X_2)$ are equal.

Using the Bers isomorphism theorem, we regard the points of the
space $\T_1$ as the pairs $X_{W^\mu} = (S_{W^\mu}, W^\mu(0))$, where
$\mu \in \Belt(\D)_1$ obey $\T_1$-equivalence (hence, also
$\T$-equivalence). Denote (for simplicity of notations) the
composition of $\mathcal J_n$ with biholomorphism $\T_1 \cong \mathcal
F(\T)$ again by $\mathcal J_n$. In view of (20) and (21), it is
presented on the fiber space $\mathcal F(\T)$ by
 \be\label{23}
\mathcal J(X_{W^\mu}) = \mathcal J(S_{W^\mu}, \ t), \quad t = W^\mu(0).
\end{equation}
This yields a logarithmically plurisubharmonic functional
$|\mathcal J_n(S_{W^\mu}, t)|$ on $\mathcal F(\T)$.

\bigskip
We have to estimate a smaller plurisubharmonic functional arising
after restriction of $\mathcal J(S_{W^\mu}, \ t)$ to $S_W \in \iota \Bigl(\fc{1}{2^p} \mathcal U^p\Bigr)$ and to $W^\mu(0)$ filling some subdomain $D_{\theta}$.

Since our functionals are polynomials, they are defined for all $S_W \in \T$ and $t$ from some domain $D_\theta$ containing $D_{\mathcal X,\theta}$.
We define on $D_\theta$ the function
$$
u_\theta(t) = \sup_{S_{W^\mu}} |\mathcal J_n(S_{W^\mu}, t)|,
$$
where the supremum is taken over all $S_{W^\mu} \in \T$ admissible for a given $t =
W^\mu(0) \in D_\theta$.

The following basic lemma is a generalization of the corresponding result in \cite{Kr5}.
It provides that this function inherits subharmonicity of $\mathcal J_n$.

\bigskip\noindent
{\bf Lemma 6}. {\it The function $u_\theta(t)$ is subharmonic on its domain $D_\theta$ filled
by the admissible values of $W^\mu(0)$. }

\bigskip
The {\bf proof} of this lemma is complicated. Similar to \cite{Kr5}, it involves the approximation
of elements from $\fc{1}{2^p} \mathcal U^p$ by polynomials given by Lemma 5 which provides the finite dimensional submanifolds wekly approximating $\iota(\big(\fc{1}{2^p} \mathcal U^p))$ in
the underlying space $\T$ (and simultaneously in the space $\T_1$) in the topology of locally uniform convergence on $\C$.

Since the set $\iota(\big(\fc{1}{2^p} \mathcal U^p))$ is a complex submanifold in $\T$, the restriction of the function $|\mathcal J(S_{W^\mu}, t)|$
to this submanifold and to the corresponding values of $t = W^\mu(0)$ also is
plurisubharmonic. The arguments from \cite{Kr5} are straightforwardly extended to this restriction, giving in a similar way  the corresponding maximal subharmonic function
$$
u_\theta(t) = \sup_{S_{W^\mu}} |\mathcal J_n(S_{W^\mu}, t)|;
$$
the supremum here is taken over $S_{W^\mu} \in \iota(\big(\fc{1}{2^p} \mathcal U^p))$.

One also has to extend the previous construction to the increasing unions of the quotient spaces
 \be\label{24}
\mathcal T_s = \bigcup_{j=1}^s \ \wh \Sigma_{\theta_j}^0/\thicksim \
= \bigcup_{j=1}^s \{(S_{W_{\theta_j}}, W_\theta^\mu(0)) \} \ \simeq
\T_1 \cup \dots \cup \T_1,
\end{equation}
where $\theta_j$ run over a dense subset $\Theta \subset [-\pi, \pi]$, the equivalence relation $\thicksim$ means $\T_1$-equivalence on a dense subset $\wh \Sigma^0(1)$ in the union $\wh \Sigma(1)$
formed by univalent functions $W_{\theta_j}(z) = e^{-i \theta_j} z +
b_0 + b_1 z^{-2} + \dots$ on $\D^*$ with quasiconformal extension to
$\hC$ satisfying $W_{\theta_j}(1) = 1$, and
$$
\mathbf W_\theta^\mu(0) := (W_{\theta_1}^{\mu_1}(0), \dots ,
W_{\theta_s}^{\mu_s}(0)).
$$
The Beltrami coefficients  $\mu_j \in \Belt(\D)_1$ are chosen here
independently. The corresponding collection $\beta = (\beta_1, \dots, \beta_s)$
of the Bers isomorphisms
$$
\beta_j: \ \{(S_{W_{\theta_j}}, W_{\theta_j}^{\mu_j}(0))\} \to
\mathcal F(\T)
$$
determines a holomorphic surjection of the space $\mathcal T_s$
onto $\mathcal F(\T)$.

Taking in each union (24) the corresponding collection $\iota_s \Bigl(\fc{1}{2^p} \mathcal U^p\Bigr)$, one obtains
in a similar fashion the increasing sequence of  maximal subharmonic functions
$$
u_s(t) = \sup_{\Theta} u_{\theta_s}(t) = \sup \big\{|\mathcal
J_n(S_{W^\mu}, t)|: \ S_{W^\mu} \in \bigcup_s \iota_s \Bigl(\fc{1}{2^p} \mathcal U^p\Bigr)\big\},
$$
whose limit
  \be\label{25}
u(t) =\lim\limits_{s\to \iy} u_s(t)
\end{equation}
is determined  and subharmonic on a disk
 \be\label{26}
D_\rho = \bigcup_{\Theta} D_{\rho,\theta_s},
\end{equation}
because the union of spaces (15) admits the circular symmetry.

\bigskip\noindent
{\bf Step 4: Determination of the range domain of $W^\mu(0)$}.
Our goal now is to find the domain of admissible values of $W^\mu(0)$, i.e. the radius
of the disk (26).
This requires a covering estimate of Koebe's type given by the following lemma.

Let $G$ be a domain in a complex Banach space $X = \{\mathbf x\}$
and $\chi$ be a holomorphic map from $G$ into the universal
Teichm\"{u}ller space $\T$ modeled as a bounded subdomain of $\B$.
Consider in the unit disk the corresponding Schwarzian differential
equations
 \be\label{27}
S_w(z) = \chi(\x)
\end{equation}
and pick their univalent solutions $w(z)$ satisfying $w(0) =
w^\prime(0) - 1 = 0$ (hence $w(z) = z  + \sum_2^\infty a_n z^n$).
Set
 \be\label{28}
|a_2^0| = \sup \{ |a_2|: \ S_w \in \chi(G)\},
\end{equation}
and let
$$
w_0(z) = z + a_2^0 z^2 + \dots
$$
be one of the maximizing functions for $a_2$.

\bigskip\noindent
{\bf Lemma 7}. \cite{Kr4} {\it (a) For every indicated solution $w(z) = z + a_2 + \dots$ of the differential equation (17), the image domain $w(\D)$ covers entirely the disk $\{|w| < 1/(2 |a_2^0|)\}$.

The radius value $1/(2 |a_2^0|)$ is sharp for this collection of functions, and the circle $\{|w| = 1/(2 |a_2^0|)$ contains points not belonging to $w(\D)$ if and only if $|a_2| = |a_2^0|$ (i.e., when $w$ is one of the maximizing functions).

(b) The inverted functions
$$
W(\zeta) = 1/w(1/\zeta) = \zeta - a_2^0 + b_1 \zeta^{-1} + b_2 \zeta^{-2} + \dots
$$
map the disk $\D^*$ onto a domain whose boundary is entirely contained in the disk} $\{|W + a_2^0| \le |a_2^0|\}$.

\bigskip
Now we show that in the case of nonvanishing $H^p$ functions this radius $2 |a_2^0|$ is naturally
connected with the extremal function $\kappa_{1,p}(z)$  maximizing the coefficient $|c_1|$.

Consider the collection $\mathcal N_p \ (p > 1)$ of all nonvanishing $H^p$ functions located
in the ball $\{\|\vp\| < 2\}$  in $\B$ and denote the minimal radius of the balls in $H^p$
containing these functions by $r(p)$; that is
$$
r(p) = \sup \{\|f\|_p: \ \|f\|_\B \le 2, \ f(z) \ne 0 \ \text{in} \ \D\}.
$$
For any such $f$, the solutions $w(z)$ of the equation $S_w = f$ are univalent holomorphic functions
on the disk $\D$. The set $\fc{1}{2^{1/p}} \mathcal U^p$ applied earlier is a proper subset
of $\mathcal N_p$.

\bigskip\noindent
{\bf Lemma 8}. {\it For any space $H^p, \ p > 1$, and its subset $\mathcal N_p$, we have the equality
  \be\label{29}
S_{w_0}(z) = r(p) \kappa_{1,p}(z)
\end{equation}
which means that the Schwarzian of the extremal univalent function $w_0(z)$ maximizing the second
coefficient $a_2$
on the set $\mathcal N_p$ equals the extremal function for $c_1$ (hence, the maximizing function
for (28) also is unique). }

\bigskip\noindent
{\bf Proof}. In view of Lemma 2, it is enough to establish that
 \be\label{30}
S_{w_0}^\prime(0) = c_1^0 \ne 0
\end{equation}
(in other words, that the zero set of the functional
$J_1(f) = c_1$ is separated from the set of rotations (9) of the function $w_0$). This yields that the maximal function (25) for the functional $|J_1(f)| = |c_1|$ is defined on the whole disk $\D_{2|a_2^0|}$, attaining its maximum on the boundary circle.

We pass to intersections
$$
B_{1,M}^0(H^p) = B_1^0(H^p) \cap \{f \in L_\iy(\D): \ \|f\|_\infty < M\},
$$
with $M < \iy$, getting the corresponding subharmonic functions
$$
u_M(t) = \sup \{|\mathcal J(S_{W^\mu}, t)|: \
S_{W^\mu} \in \iota \Bigl(\fc{1}{2^p} \mathcal U^p \cap \mathcal P \Bigr) \cap B_{1,M}^0(H^p)\}
$$
with $\lim\limits_{M\to \infty} u_M(t) = u(t)$ and the  points $f_M = S_{w_{0,M}}$, maximizing $|a_2|$
on these sets. The collection $\{f_M\}$ is weakly compact in $H^p$.

Applying Lemma 7, one obtains similar to \cite{Kr5} that both maximal values $|a_2^0|$ and
$|c_1^0|$ are obtained on the same function
$$
f_0(z) = \lim\limits_{M \to \infty} f_M(z) = S_{w_0},
$$
and the uniqueness in  Lemma 2 yields that this function must coincide with $r(p) \kappa_{1,p}$.
This completes the proof of Lemma 8.

\bigskip\bigskip
{\bf Step 5: Finishing the proof}. Now we can prove the assertion of the theorem.
The assumption $p \ge 2$ insures that the boundary function
$f(e^{i \theta}) = \lim_{r\to 1} f(r e^{i \theta})$ of any $f \in H^p$ admits Parseval's equality
  \be\label{31}
1 \ge \frac{1}{2 \pi} \int\limits_{- \pi}^\pi |f(e^{i \theta})|^2 d \theta
= \sum_1^\infty |c_n|^2.
\end{equation}
In particular, for $f(z) = \kappa_{1,p}(z) = \sum\limits_0^\iy c_n^0 z^n$
we have from (2)
 \be\label{32}
|c_1^0|^2 = (2/e)^{2(1-1/p)} = 0.5041...^{1-1/p} > 0.5041...
\end{equation}
for all $p > 1$. Hence,  by (31),
 \be\label{33}
\sum_2^\infty |c_n^0|^2 < 0.5 < |c_1^0|^2.
\end{equation}

\bigskip
Now take $n = 2$ and, letting $f_2(z) = f(z^2)$, consider on the set $B_1^0(H^p)$ the  functional
$$
I_2(f) = \max \ (|J_2(f)|, |J_2(f_2)|).
$$
Since the correspondence $f(z) \mapsto f_2(z)$ is linear, the functional $J_2(f_2)$ is holomorphic
with respect to
$f$ in $H^{2m}$ norm and naturally extends to a holomorphic functional on the spaces $\T$.
Hence, the functional $I_2$ is plurisubharmonic on $\T$. It lifted to the covering space $\T_1$
together with $J_2$.

Similar to above, this lifting generates via (25) a nonconstant radial subharmonic function $u_2(t)$
on the disk $\{|t| < 1/2|a_2^0|\}, \ t = W^\mu(0)$.
This function is logarithmically convex, hence monotone increasing, and thus attains its maximal
value at $|t| = 2 |a_2^0|$.

Taking into account the connection between the extremal value $|a_2^0|$ and the function $\kappa_{1,p}$
established by Lemma 8, one concludes  that the maximal value of $I_2(f)$ on $B_1^0(H^p)$ is attained
on the pair $(f, f_2)$ with
$$
f(z) = \kappa_{1,p}(z), \quad f_2(z) = \kappa_{1,p}(z^2).
$$
Since the set of admissible maps $w(z) \in \wh S(1)$ with $S_w = f$ for $J_2(f_2)$ is the same as for $J_2(f)$, one derives from above
$$
\max_{B_1^0(H^p)} I_2(f) = \max \ \{|c_1^0|, |c_2^0|\},
$$
which by (33) is equal to $= |c_1^0| = (2/e)^{1-1/p}$. This yields the desired estimate (1) for $n = 2$;
the extremal maximizing function is determined up to the pre and post rotations about the origin.

Now consider subsequently for $n = 3, 4, \dots$ the functionals
$$
I_n(f) = \max \ \{|J_n(f)|, |J_n(f_2)|,  \dots, |J_n(f_n)|\}
$$
with $f_n(z) = f(z^n)$. Similar to $I_2$, this does not expand the set of admissible maps $w(z) \in \wh S(1)$
with $S_w = f$, and therefore, $|I_n(f)|$ has the same maximum, as $|J_n(f_n)|$.

Each functional $I_n$ generates similar to above the corresponding circularly symmetric subharmonic
function $u_n(t)$ on the disk $\{|t| < 1/2|a_2^0|\}, \ t = W^\mu(0)$, which provides in the same
way the bound
$$
\max_{B_1^0(H^p)} I_n(f) = \max \ \{|c_1^0|, |c_n^0|\} = (2/e)^{1-1/p},
$$
with a similar description of the extremal functions.
This completes the proof of Theorem 1.

\bigskip\bigskip
\centerline{\bf 4. ADDITIONAL REMARKS}

\bigskip\noindent
{\bf 4.1. Holomorphy in parameters}.
As was mentioned in the proof of Theorem 1, the holomorphic dependence of normalized quasiconformal
maps on complex parameters is an underlying fact for the Teichm\"{u}ller space theory and for many
other applications. It was first established and applied by Ahlfors and Bers in \cite{AB} for maps with
three fixed points on $\hC$.

Another somewhat equivalent proof of holomorphy involves the variational
technique for quasiconformal maps. For the maps $w$ from $S_{\theta,\a}$,
this holomorphy is a consequence of the following lemma  from \cite{Kr1}, Ch. 5
combined with appropriate M\"{o}bius maps.

\bigskip\noindent
{\bf Lemma 9}. {\it Let $w(z)$ be a  quasiconformal map of the plane
$\hC$ with Beltrami coefficient $\mu(z)$ which satisfies
$\|\mu\|_\iy < \ve_0 < 1$ and vanishes in the disk $\{|z| < r\}$.
Suppose that $w(0) = 0, \ w^\prime(0) = 1$, and $w(1) = 1$. Then,
for sufficiently small $\ve_0$ and for $|z| \le R < r_0(\ve_0, r)$
we have the variational formula
$$
w(z) = z - \frac{z^2 (z - 1)}{\pi} \iint\limits_{|\zeta|>r} \
\frac{\mu(\zeta) d\xi d \eta}{\zeta^2(\zeta -1)(\zeta- z)} +\Omega_\mu(z),
$$
where $\zeta = \xi + i \eta; \ \max_{|z|\le R} |\Omega_\mu(z) \le C(\ve_0, r, R) \|\mu\|_\iy^2; \
r_0(\ve_0, r)$ is a well defined function of $\ve_0$ and $r$ such that
$\lim_{\ve_0\to 0} \ r_0(\ve_0, r) = \infty$, and the constant $C(\ve_0, r, R)$
depends only on $\ve_0, \ r$ and $R$.   }

\bigskip\noindent
{\bf 4.2 Remarks on the case $1 < p < 2$}.
All arguments in the proof of Theorem 1, excluding the Parseval equality, applied in the last step,  work for any $p > 1$.
In fact, this equality was applied only to the function $\kappa_{1,p}$ maximizing
$|c_1|$ and was used by estimation $|c_n|$ for comparison of the initial non-free coefficient of functions $\kappa_{1,p}(z^m), \ 1 \le m \le n$.

The explicit representation (2) of $\kappa_{1,p}$ shows that this function is bounded on the unit disk if $1 < p < 2$;
hence belongs to $H^2$. However, since for $1 < p < 2$,
$$
\|\kappa_{1,p}\|_{H^2} > \|\kappa_{1,p}\|_{H^p}
$$
the needed relation (31) giving (32), (33) fails. The functionals $J_n$ and $I_n$ cannot be compared on this way.

\bigskip\noindent
{\bf 4.3. On extremal functions in Bergman spaces}.
One of the interesting extensions of the  Hummel-Scheinberg-Zalcman problem
(also still unsolved) is to estimate the Taylor coefficients of
nonvanishing holomorphic maps $f(z) = c_0 + c_1 z + \dots$ of the unit disk $\D$ into other complex Banach spaces $X$.
Denote by $B(X)$ the unit ball of $X$.

We illustrate here on the case of Bergman's space $A_2$ that the features of extremal functions can be essentially different from above. Recall that the norm of $A_2$ is  $\|f\| = (\fc{1}{\pi} \iint_\D |f(z)|^2 dx dy)^{1/2}$.

The collection $B_0(A_2)$ of nonvanishing holomorphic functions $f(z)$  mapping the
disk $\D$ into the closed ball $\ov{B(A_2)}$ (i.e., with $f(z) \ne 0$ on $\D$ and $\|f\| \le 1$) is compact
in the weak topology of the locally uniform convergence in $\D$. So any holomorphic coefficient functional
 \be\label{34}
J(f) = J(c_{m_1}, \dots, c_{m_s}) \quad \text{with} \ \ 1 \le m_1 < m_2 < \dots < m_s = N < \iy
\end{equation}
has an extremal $f_0$ on which $|J(f)|$ attains its maximum on $B_0(A_2)$.

While the extremal functions of many problems in Hardy spaces are bounded,
Proposition 2 implies, for example, that {\it any function $f_0 \in B_0(A_2)$ maximizing the functional (34)
must be unbounded on the disk $\D$} (compare with the extremal problems
for nonvanishing Bergman functions investigated e.g. in \cite{ABKS}, \cite{BK}).

This difference is caused by the fact mentioned after the proof of Proposition 2: quasiconformal deformations
created by Propositions 1 and 2 preserve the norm in $A_p$, while the norm of the Hardy spaces can be increased.

Indeed, it follows from Proposition 2 that any extremal $f_0$ of $J(f)$ on $B_0(A_2)$  must be unbounded on $\D$, unless $f_0$ is a zero-free polynomial
 \be\label{35}
p_N(z) = c_0 + c_1 z + \dots + c_N z^N
\end{equation}
(with $c_0 \ne 0$); otherwise, one can vary the coefficients $c_k$ and obtain by this lemma an admissible function $f_{*} \in B_0(A_2)$ with $|J(f_{*})| > |J(f_0)|$.

It remains to establish that the polynomials (35) with $\|p_N\|_{A_2} \le 1$ cannot be extremal for $J(f)$.
We pick a sufficiently small $\ve > 0$ and consider
the polynomial
$$
p_\ve(z) = -\ve c_0 + \ve z^{N+1},
$$
for which
$$
\max_{\mathbb S^1} |p_\ve(z)| < \max_{\mathbb S^1} |p_N(z)|.
$$
Then the Rouch\'{e} theorem yields that the polynomial
$$
P_{N+1,\ve}(z) = p_N(z) + p_\ve(z) = (1 - \ve) c_0 + c_1 z + \dots + c_N z^N + \ve z^{N+1}
$$
also must be, together with $p_N$, zero-free on $\D$. Its norm is estimated by
$$
\|p_{N+1,\ve}\|_{A_2}^2 = (1 - \ve)^2 |c_0|^2 + \fc{|c_1|^2}{2} + \dots
+ \fc{|c_N|^2}{N + 1} + \fc{\ve^2}{N + 2} = \|p_N\|_{A_2}^2 - 2 \ve + O(\ve^2) <   \|p_N\|_{A_2}^2 = 1,
$$
which implies that $P_{N+1,\ve}$ is an admissible function, with
 \be\label{36}
|J(p_{N+1,\ve})| = |J(p_N)| = \max \{|J(f)|: \ f \in B_0(A_2)\}.
\end{equation}
But this contradicts to Proposition 2, because this proposition allows one to construct the variations of $p_{N+1,\ve}$, which preserve its $A_2$-norm and increase $|J(p_{N+1,\ve})|$, disturbing (36).
This completes the proof of our claim.

\medskip
{\small\em{ \leftline{Department of Mathematics, Bar-Ilan
University, Ramat-Gan 5290002, Israel} \leftline {Department of Mathematics,
University of Virginia, Charlottesville, VA 22904-4137, USA}}}

\end{document}